\documentclass[11pt]{amsart}
\usepackage{amssymb}
\usepackage{amsmath}
\newtheorem{lemma}{Lemma}
\newtheorem{theorem}{Theorem}
\newtheorem{proposition}{Proposition}
\newtheorem{corollary}{Corollary}

\theoremstyle{definition}

\newtheorem{example}{Example}

\newtheorem{remark}{Remark}

\def\det{{\rm{det}}}

\newcommand{\Hom}{\mbox{Hom}}

\newcommand{\ignore}[1]{}
\newcommand{\can}{\overline{\phantom{x}}}

\keywords{Flexible quadratic algebras, composition algebras, cross products.}

\subjclass[2000]{Primary: 17A45; Secondary: 17A20, 17A75, 11E25}

\title{On flexible quadratic algebras}

\author{S. Pumpl\"un}

\address{School of Mathematics\\
 University of Nottingham\\
 University Park\\
 Nottingham NG7 2RD\\
 United Kingdom\\
 susanne.pumpluen@nottingham.ac.uk
 }

\begin{document}

\maketitle

\begin{abstract} Let $R$ be a ring. A construction method for flexible quadratic algebras
with scalar involution over $R$ is presented  which unifies various classical
constructions in the literature, in particular those to construct composition algebras.
\end{abstract}


\section*{Introduction}

Let $R$ be a unital commutative associative ring.
We present a construction method for $R$-algebras starting from an associative quadratic $R$-algebra with a scalar involution.
 Special cases of this construction yield flexible quadratic algebras over rings
 (an algebra $A$ is {\it flexible} if $x(yx)=(xy)x$ for all elements $x,y\in A$).
 Composition algebras are certainly among the best-known classes of flexible
quadratic algebras. Our method unifies the different approaches by Petersson [P2], Thakur [T] and [Pu]
to construct composition algebras. In special cases  it overlaps
with the one described by Becker [B] which yields flexible quadratic division algebras
 containing subalgebras of half their rank.

\section{Preliminaries}

 Let $A$ be an $R$-algebra.
The term ``$R$-algebra" always refers to a unital nonassociative
algebra $A$ which is faithfully projective as an $R$-module [K, p.~52].

 $A$ is called {\it alternative} if its associator $[x, y, z]= (xy)z - x(yz)$ is
alternating. The {\it nucleus} of $A$ is
defined as $N(A)  = \{ x \in A \, \vert \, [x, A, A] = [A, x, A] = [A,
A, x] = 0 \}$. The nucleus is an associative subalgebra of $A$ (it may be
zero), and $x(yz) = (xy) z$ whenever one of the elements $x, y, z$ is in $N(A)$.
 An anti-automorphism $\sigma:A\to A$ of period 2 is called an {\it involution} on $A$.
If $2$ is an invertible element in $R$, we have $A = {\rm Sym}(A,\sigma) \oplus
{\rm Skew}(A,\sigma)$ with
${\rm Skew}(A, \sigma)= \{ x \in A \;\vert \; \sigma (x) = - x \}$
the set of skew symmetric elements and
${\rm Sym}(A, \sigma) = \{ x \in A \; \vert \; \sigma (x) = x \}$
the set of symmetric elements in $A$ with respect to $\sigma$.
An involution is called {\it scalar} if all {\it norms} $\sigma(x)x$ are elements of $R1$.
For every scalar involution $\sigma$, $n_A(x)=\sigma(x)x$ (resp. the {\it trace} $t_A(x)=\sigma(x)+x$) is a quadratic (resp. a
linear) form on $A$. $A$ is called {\it quadratic}, if there exists a quadratic
form $n \colon A \to R$ such that $n(1_A) = 1$ and $x^2 - n(1_A, x)x + n(x)
1_A = 0$ for all $x \in A$, where $n(x,y)$ denotes the
induced symmetric bilinear form $n(x, y) = n (x+y) - n(x) -n(y)$. The form $n$ is uniquely determined
and called the {\it norm} of the quadratic algebra $A$ [P2].
 The existence of a scalar involution  on an algebra $A$ implies that $A$ is quadratic [M1].

An $R$-algebra $C$ is called a
{\it composition algebra} if it carries a quadratic form $n \colon C \to R$
whose induced symmetric bilinear form $n(x, y)$ is {\it nondegenerate} (i.e., determines an $R$-module isomorphism $C
\stackrel{\sim}{\longrightarrow} C^{\vee} = \Hom_R (C, R)$) and which satisfies
 $n(xy) = n(x) n(y)$ for all $x, y \in C$ (it {\it permits composition}).
 Composition algebras are quadratic alternative; any nondegenerate quadratic form $n$ on the composition algebra
 which permits composition is its norm as a quadratic algebra and thus is unique [P2]. It
is called the {\it norm} of $C$ and is also
denoted by $n_C$. A quadratic alternative algebra is a composition algebra if
and only if its norm is nondegenerate [M1, 4.6]. Composition algebras only
exist in ranks 1, 2, 4 or 8. Those of rank~2 are exactly the
quadratic \'etale $R$-algebras, those of rank 4 exactly the well-known
 quaternion algebras. The ones of rank 8 are called {\it octonion algebras}. A composition algebra $C$
has a {\it canonical involution} $\can$ given by
$\overline{x} = t_C(x)1_C - x$, where $t_C \colon C \to R$, $t_C(x) = n (1_C, x)$,
is the {\it trace} of $C$. This involution is scalar.

\section{Flexible quadratic algebras}

Let $N$ be a finitely generated
projective $R$-module of constant rank. An alternating $R$-bilinear map
$\times: N\times N \to N$ is called a {\it cross product} on $N$.
 Let $D$ be an associative $R$-algebra with a scalar involution $\sigma=\can$ (e.g.
 a quadratic \'etale or a quaternion algebra over $R$ with canonical involution).
 Let $F$ be a locally free right $D$-module of constant finite
 rank $s$ together with a sesquilinear form
 $h:F\times F\to D$ (i.e., $h$ is a biadditive map such that $h(ua,vb)=\bar ah(u,v)b$ for all $a,b\in D$, $u,v\in F$). Let $\times$ be a cross product on $F$, where now $F$
 is viewed as an $R$-module, i.e. together with
$\times: F\times F \to F$, $F$ is an anticommutative $R$-algebra.
 The $R$-module $A=D\oplus F$ becomes a unital $R$-algebra denoted by $A=(D,F,h,\times)$ via the multiplication
$$(\alpha,u)(\beta,v)=(\alpha\beta-h(v,u), v \alpha +u \overline{\beta}
 + v\times u),$$
 for all $\alpha, \beta\in D$, $u,v\in F$. The algebra $D$ is a subalgebra of $(D,F,h,\times)$.

Define $$\begin{array}{l}\sigma_A:A\to A, \quad (\alpha,u)\to (\overline{\alpha}, -u),\\
n_A:A\to D, \quad n_A((\alpha,u))= \sigma_A(\alpha,u)(\alpha, u),\\
t_A:A\to R, \quad  t_A((\alpha,u))= \sigma_A(\alpha,u)+(\alpha, u)=
(t_{D}(\alpha),0).
\end{array}$$
Obviously, $\sigma_A$ is an $R$-module isomorphism of order 2. We have
$$\begin{array}{l} n_A((\alpha,u))= n_D(\alpha)+h(u,u),\\
 n_A((\alpha,u))=
(\alpha,u) \sigma_A(\alpha,u),\\
 n_A((\sigma(\alpha,u))= n_A((\alpha,u))
\,\, {\rm and}\,\, t_A(\sigma(\alpha,u))= t_A((\alpha,u)).
\end{array}$$
The map $n_A$ is quadratic with associated $R$-bilinear map
$n_A:A\times A
\to D$ given by
$$\begin{array}{l}
 n_A((\alpha,u),(\beta,v))=n_D(\alpha,\beta)+(h(u,v)+h(v,u))=\\
 (\alpha,u)\sigma_A(\beta,v)+(\beta,v)\sigma_A(\alpha,u).
 \end{array}$$
 Obviously, $F={\rm ker}(t_A)={\rm ker}(t_D)\oplus F$ and $u\times v=uv-\frac{1}{2}n_A(u,v)$.

\begin{lemma}
(i) $\sigma_A$ is a scalar involution  if and only if $h$ is a hermitian form (i.e., $h(u,v)=\overline{ h(v,u)}$
for all $u,v\in F$).\\
(ii) If $h$ is a hermitian form, then
 $t:A\times A\to R,\, t(x,y)=t_A(xy)$ is a symmetric $R$-bilinear form.
\end{lemma}

Lemma 1 is a weak  generalization of [B, Lemma 1].

\begin{corollary} If $h:F\times F\to D$ is a (perhaps degenerate) hermitian form, then
  $A=(D,F,h,\times)$ is a quadratic $R$-algebra with scalar involution
$\sigma_A:A\to A, \,\sigma_A((\alpha,u))=(\overline{\alpha},-u)$
and norm $n_A:A \to R,x\to x\sigma_A(x)$, where $n_A((\alpha,u))=n_D(\alpha)+h(u,u)$.
Moreover, if  $n_A$ is isotropic (e.g., if $h$ is isotropic) then $A$ has zero divisors.
 \end{corollary}

Note that even if $D$ is an associative composition division algebra over a field, the fact that
two quadratic algebras $A=(D,F,h,\times)$ and $A'=(D,F',h',\times')$ have isometric norms $n_A\cong n_{A'}$
only implies that the nondegenerate hermitian forms $(F,h)$ and $(F',h')$ are isometric and yields no
information on the cross products
$\times$ and $\times'$.

 Quadratic $R$-algebras with
 a scalar involution were constructed using finite dimensional vector spaces
  which carry a bilinear form $B$ by Osborn [O]. Finite-dimensional quadratic division algebras over a
 field of characteristic not 2 were shown to be
 even-dimensional by Petersson [P1]. Already earlier, Kuz'min [Ku] proved that
any quadratic division algebra, over a field of characteristic not 2 in which any
two elements generate a subalgebra of dimension at most four, must have finite dimension
and that, indeed, only the dimensions 1, 2, 4 or 8 are possible. Independently,
Kunze and Scheinberg [KS, Theorem 2.4] enumerated all unital alternative algebras over fields of
characteristic not 2 - including the infinite-dimensional ones - which have a scalar involution.

\begin{remark}
 For $A=(D,F,h,\times)$, $2\in R^\times$, we obtain ${\rm Sym}(A,\sigma_A)=R$ and
${\rm Skew}(A,\sigma_A)={\rm Skew}(D,\sigma_D) \oplus F.$
\end{remark}

\begin{lemma} If
 $h:F\times F\to D$ is a hermitian form, then $(D,F,h,\times)$ is\\
(i) flexible if and only if
 $h(u\times v,u)+\overline{h(u\times v,u)}=n_A(u\times v,u)=0$ and
$(u\times v) \times u=u\times ( v \times u)$ for all $u,v\in F$;\\
(ii) 
alternative if and only if $h(u, u\times v)=0$ and
$u\times (u\times v)=-h(u,u)v+h(v,u)u$ for all $u,v\in F$,
if and only if
$h( u\times v,v)=0$ and
$(u\times v)\times v=h(v,v)u-h(u,v)v$ for all $u,v\in F$.
\end{lemma}

If $(D,F,h,\times)$ is flexible with a scalar involution, then it is a noncommutative
 Jordan algebra [M1, (3.3)].
Due to the nature of the multiplication in
$(D,F,h,\times)$ it suffices to work with elements of $F$
to check flexibility. The identities in (ii) were already obtained
in [T, (2.1)] in a more restricted setting (there, F has rank 3, $D$
is quadratic \'etale, and $h$ is a nondegenerate hermitian form).
If $\times$ is the zero-map, then $(D,F,h,0)$ is trivially flexible.

\begin{remark} If $D=R$, then $\sigma=id$ and $h=B$ is a symmetric $R$-bilinear form on $F$. $A=(R,F,B,\times)$ is the
unital quadratic $R$-algebra described in [B, p.~26] (see also [O, Theorem 1]).
It is well-known that any quadratic algebra $A$ over a field of characteristic not 2 with a scalar involution can
is of the type $A=(R,F,B,\times)$ for a suitable $F$, a symmetric
bilinear form $B$ and a cross product $\times$.\\
(i) By Lemma 2 (i),  $A$ is flexible if and only if $B(u,v\times u)=B(u,u\times v)$ and
$(u\times v) \times u=u\times ( v \times u)$ for all $u,v\in F$,
which is equivalent to
$$B(u,v\times u)=0\quad {\it and}
\quad (u\times v) \times u=u\times ( v \times u)$$
 for all $u,v\in F$.
\\(ii) Let $2\in R^\times$. If we apply [B, Lemma 2, Satz 1] to the residue class algebras, we see
 that $A$ is flexible if and only if $t_A$ is associative,
if and only if $B$ is associative; and that $A$ is flexible and $t_A$ and $n_A$ are weakly nondegenerate
if and only if $B$ is  associative and weakly nondegenerate.
($q$ is {\it weakly nondegenerate} if ${\rm Rad}(q)=\{a\in A | q(a,A)=0\}=0$.)\\
 (iii) For a flexible quadratic $R$-algebra $A=(R,F,B,\times)$ of rank
greater or equal to 3 with a weakly nondegenerate associative bilinear
form $B$, the residue class algebra $A(P)$ is central-simple
for each $P\in {\rm Spec}\,(R)$. This follows from [B, Satz 1, Satz 2] applied to the residue
class algebras. If $A$ is a quadratic algebra of rank 3 over a domain $R$ and $B$  weakly
nondegenerate and associative, then $A$ is commutative. This follows from [B, p.~33].
\end{remark}



\begin{theorem} Let $2\in R^\times$.\\
(i) Let $F$ be a finitely generated projective $R$-module of constant rank $s$.
Let $B:F\times F\to R$ be a nondegenerate symmetric bilinear form
on $F$ and $\times:F\times F\to F$ a cross product, such that
\begin{equation}\tag{1}
B(u\times v,u\times v)=B(u,u)B(v,v)-B(u,v)^2
\end{equation}
\begin{equation}\tag{2}B(u,u\times v)=0 \quad
\end{equation}
for all $u,v\in F$. Then  $C=(R,F,B,\times)$
is a composition algebra over $R$.
 In particular $s=1$, $s=3$, or $s=7$, and ${\rm det}(F,N)=\bigwedge^s(F,N)=R$.\\
 (ii) Let $C$ be a composition algebra over $R$ with norm $n$ of rank $s+1\geq 1$. Then
 $C=R\oplus F$ with $F={\rm Skew}(C,\can)$.
The $R$-module $F$ satisfies ${\rm det}(F)=\bigwedge^s(F)=R$ with $s=1,3$ or 7.
 Define $\times: F\times F\to F$
via $u\times v=  \frac{1}{2}[u,v]= \frac{1}{2} (uv-vu)$
and $B:F\times F\to R$ via $B(u,v)=-\frac{1}{2}(uv+vu)=
\frac{1}{2}(u\bar v+v\bar u)=\frac{1}{2}n(u,v)$.
Then $\times$ is a cross product and $B$ is a nondegenerate
symmetric bilinear form on $F$ satisfying equations $(1)$ and $(2)$ in (i) for all $u,v\in F$.
\end{theorem}

This is well-known if $R$ is a field. Lemma 2 (ii) implies the first part of Theorem 1. Otherwise, the
proof of [Mo, Theorem 10] can be used verbatim.
Note that the intrinsic value of this theorem is impaired by the fact that producing a cross product on a finitely generated
projective $R$-module $F$ of rank 7 satisfying (1) and (2) is probably not any easier than producing an octonion algebra
structure on $R\oplus F$.

If, in (ii), $C$ is quadratic \'etale, then there exists a line bundle $L$ over $R$ and a nondegenerate quadratic form
$N_L:L\times L\to R$ such that $C\cong{\rm Cay}(R,L,N_L)$. Since every bilinear map $L\times L\to L$ is symmetric,
this forces the cross product on $L$ to be zero.

\begin{example}  \ignore{
(due to H.P. Petersson) Let $k$ be an infinite field of characteristic not 2.
Suppose $A$ is a flexible quadratic algebra over $k$ such that the nil radical of $A$ does not have codimension 1.
Then $A$ contains a quadratic  \'etale subalgebra $D=k[x]$, but $A\not\cong(D,F,h,\times)$ for suitable $D$,
$F$, $h$ and $\times$:
otherwise - by a simple calculation involving the definition of the multiplication of $A$ as described on p.~3
- $A$ would be left-alternative, hence alternative,
a contradiction.}
 Let $A$ be a quadratic $R$-algebra  of rank 2. Then $A$ is commutative associative and the conjugation is a
scalar involution (see [K, I.(1.3.6)] or check it locally). We get $A=(R,L,B,\times)$ for a suitable $L\in {\rm Pic}\,R$,  $B:
L\times L\to R$ a bilinear form. $A$ is commutative, therefore $\times$ is
the zero-map, and $A={\rm Cay}(R,L,N)$ with $N:L\to R$ the - perhaps degenerate -
quadratic form associated to $B$.
 \end{example}

\begin{lemma} Let $2\in R^\times$. Let $A$ be a flexible quadratic $R$-algebra with
a scalar involution $\sigma$ which contains a composition subalgebra $D$
of rank 2 and no other associative composition subalgebra of larger rank, or
a composition subalgebra $D$ of rank 4. Let $F=D^\perp$.
 Define $-h$ to be the orthogonal projection of the multiplication $vu$ for $u,v\in F$ onto $D$ and $\times$
  to be the orthogonal projection of the multiplication $vu$ for $u,v\in F$ onto $F$. Then:
\\ (i) $\times:F\times F\to F$ is a cross product.
\\ (ii)  $n_A(u,v)=t_{D}(h(u,v))$ and $h(u,v)=\overline{h(v,u)}$.
\ignore{
 (iii) If the action $F\times D\to F$,
$(\alpha,u)\to u\alpha$
given by multiplication in $A$ makes $F$ into a locally free right $D$-module,
 $\times$ is sesquilinear and $(ub)v=(uv)b$ for all $u,v\in F$, $b\in D$,
then $h:F\times F\to D$ is a hermitian form and $A\cong(D,F,h,\times)$.
}
\end{lemma}

 \begin{proof} Let $A$ be a quadratic $R$-algebra with a scalar involution $\sigma$, which contains
 a composition subalgebra $D$ of rank 2 or 4 as subalgebra of maximal rank with nondegenerate norm.
  Since the norm $n_A(x)=x\sigma(x)$
 is nondegenerate when restricted to $D$, we can decompose $A$ as $A=D\perp D^{\perp}$ [K, (3.6.2), p.~17]. Since
  $-h$ and $\times$ are the orthogonal projections
 of the multiplication $vu$ for $u,v\in F$ onto $D$ and onto $F$, we have $vu=-h(u,v)+u\times v$ for $u,v\in F$.
\\ (i) Since $F\subset {\rm Skew}(A,\sigma)$, it follows that $u^2\in R$ for all
$u\in F$ and thus $u\times u=u^2+h(u,u)\in D \cap D^{\perp}=\{0\}$. Hence $\times$ is
anticommutative and a cross product on $F$.
\\ (ii) We have $n_A(u)=-u^2=h(u,u)$ for all $u\in F$. Since
$u\times v=-v\times u$, it follows that $h(u,v)+\overline{h(u,v)}=-vu-uv=v\sigma(u)+
u\sigma(v)=n_A(u,v)$, i.e. $n_A(u,v)=t_{D}(h(u,v))$ and $h(u,v)=\overline{h(v,u)}$.
\ignore{
 (iii) If the action $F\times D\to F$, $(\alpha,u)\to u\alpha $
given by multiplication in $A$ makes
$F$ into a locally free right $D$-module, then $h$ is a hermitian form if and only if
$h(u,vb)=h(u,v)b$ for all $u,v\in F$, $b\in D$, which is equivalent to $-(vb)u+u\times (vb)=
-(vu)b+(u\times v)b$. Hence, if $\times$ is sesquilinear as well as $(ub)v=(uv)b$ for all $u,v\in F$,
 $b\in D$, then $h$ is hermitian.

--------------------------------------------

WHY should this end up in $F$, not $A$?
}
\end{proof}


\begin{example}
(a) Let $2\in R^\times$. Let $A=(D,F,h,\times)$ be an algebra with a scalar involution,
 $D$  quadratic \'etale or a quaternion algebra. The algebra $\hat{A}=(R,F,B,\times)$,
 $B(u,v)=\frac{1}{2}(h(u,v)+\overline{h(u,v)})=\frac{1}{2}
n_A(u,v)$, is flexible if and only if $A$ is flexible.
If $A$ is alternative, the identity $((u\times v)\times v)\times v=
n_A(v) u\times v$ holds in $\hat{A}$ for all $u,v\in F$ (Lemma 2).
In particular, if $A={\rm Cay}(T,F,h,\times_{\alpha})$ is an octonion algebra over
 a field $R$ of characteristic not $2$, then
$\hat{A}=(R,F,B,\times_{\alpha})$ is the colour algebra defined in [E].
\\(b) If the cross product is trivial,  $A=(D,F,h,0)$ with
$S$ a quadratic  \'etale $R$-algebra or a quaternion algebra over $R$, is
a unital noncommutative Jordan algebra. In particular,
$(R,F,-B,0)={\rm JSpin}(F,B)$ is the well-known algebra of the Spin Factors
 [M2, p.~178].\\
(c) Let $k$ be a field of characteristic not 2.  Given a quadratic $k$-algebra
$A$ with scalar involution $*$, the map $(a,b)^*=(a^*,-b)$ defines a scalar
involution  on the quadratic algebra $C(A,\cdot_1,\cdot_2,\cdot_3)
=A\oplus A$ with multiplication $(a,b)(c,d)=(ac+d^* \cdot_1 b,d\cdot_2 a+b\cdot_2 c^*
 +b\cdot_3 d)$, where the maps $\cdot_i$ are certain bilinear products on $A$ [B].
 Using this construction all flexible quadratic division algebras of dimension $2m$ over $k$
 which contain a unital subalgebra of dimension $m$
 can be obtained [B, p.~41]. The construction generalizes
to quadratic algebras with scalar involutions over arbitrary rings where
2 is invertible.  However,  it will not produce all such algebras.
\ignore{, since there exist
octonion (or even quaternion) algebras over rings, which do not contain any
 subalgebras with scalar involutions other that the base ring itself.}
If a flexible quadratic algebra over a field contains a composition
subalgebra of dimension 2 or 4, the method coincides
with ours whenever the algebra has twice its dimension, i.e. dimension
4, respectively 8. In that case $ d^* \cdot_1 b$ is given by the hermitian
form $h(d,b)$, $\cdot_2$ is the usual algebra multiplication of $A$,
and $ \cdot_3$ is our cross product $\times$.
\end{example}

Let $S$ be  a quadratic  \'etale $R$-algebra. For the remainder of this section, we restrict our attention to the
special case that $F$ is a finitely generated projective right $S$-module of constant rank  carrying a
nondegenerate hermitian form $h \colon F \times F \to S$.
 For each invertible $\mu \in S$, the form  $\mu h \colon F \times F \to S$
 is a nondegenerate $\epsilon$-hermitian form, with $\epsilon=
 \overline{\mu}/
 {\mu}$. In particular, $\epsilon\not =1$ if and only if $\mu \not\in R$,
 and $\epsilon=-1$ if and only if $\mu \in {\rm Skew}(S,{\can})$.

Recall that conversely, the study of $\varepsilon$-hermitian forms over $S$
can be reduced to the study of hermitian ones by scaling in several important
cases:
Let $\mu \in S^{\times}$ such that $\frac{\overline{\mu}}{\mu} = \varepsilon$.
Then $\mu h$ is a hermitian form, for any $\varepsilon$-hermitian form $h$. Such
a $\mu$ exists if $H^1 ({\Bbb Z}/2 {\Bbb Z}, S^{\times}) = 0$ (Hilbert's
Theorem 90), e.g. if ${\rm Pic}\, R = 0$ [K, p.~300].

\begin{proposition} Let $A=(S,F,\mu h,\times)$ for some nondegenerate hermitian form $h:F\times F\to S$, and an
invertible scalar $\mu\in S\setminus R$.
 If $S$ is a domain (or if $F$ is torsion-free and $\bar\mu-\mu$ not a zero divisor in $S$), then $N(A)=S$.
\end{proposition}

\begin{proof}
A straightforward calculation shows that $S \subset N(A)$.

 Now let $(e, e') \in N(A)$ with $e \in S$ and $e' \in F$.
We have to show that this implies $e' = 0$. The equation
 $$(e, e')((u, w) (u', w')) = ((e, e') (u, w)) (u', w')$$
 implies that
$$e' \bar\mu h( w,w') = w'  \mu h (w, e')$$
 for all $w, w' \in F$ and thus (put $w'=e'$)
$$e' (\bar\mu-\mu)h(w,e')=0$$
 for all $w \in F$. If $h(w,e')=0$ for all $w\in F$ then $e'=0$ since $h$ is nondegenerate
and thus $(e,e')=(e,0)\in S$. Otherwise, there is one $w\in F$ such that $h(w,e')\not=0$. Then
if $S$ has no zero divisors, i.e., is a domain, $F$ is torsion free and we have $e'=0$ since $\bar\mu-\mu\not=0$ here.
The same holds if $F$ is torsion free and $\bar\mu-\mu$ not a zero divisor in $S$.
\end{proof}

\section{Constructions of composition algebras}

 Petersson [P2] introduced the following {\it generalized Cayley-Dickson
 doubling process} for composition algebras over rings: Let $D$ be a composition algebra of rank $\leq 4$ over $R$ with
 canonical involution $\can$.
Let $P$ be a locally free right $D$-module of rank one and norm one (cf. [P2, 2.3] for the definition; if
$D$ is a quaternion algebra, of norm one means that the reduced norm functor defined in [KOS] associates to $P$
 the free line bundle $R$).
 Thus there exists a nondegenerate quadratic form $N \colon P \to R$
satisfying $N(w \cdot u)=N(w)N_D(u)$ for all elements $w\in P$,
$u\in D$, where $\cdot$ denotes the right $D$-module
structure of $P$ (for $D$ a quaternion algebra, see also [K, p.~161 ff]). $N$ is uniquely determined up to an invertible
factor in $R^\times$ and determines a unique $R$-bilinear map $P \times P \to D$, written multiplicatively and satisfying
 $(w \cdot u)(w \cdot v)=N(w)\overline{ v}u$ for all $u,v\in D$, $w\in P$.
The $R$-module ${\rm Cay} (D, P, N) = D \oplus P$
becomes a composition algebra under the multiplication
$$(u, w) (u', w') = (u u' +  ww', w' \cdot u + w \cdot \overline{u}')$$
for $u, u' \in D$, $w, w' \in P$, with $\cdot$ denoting the right $D$-module structure
of $P$.  Its norm is given by
$$n_{{\rm Cay} (D, P, N)}((u,w))=n_D(u)-N(w).$$
If $F$ a
locally free right $D$-module of rank 1 and $h:F\times F\to D$ a nondegenerate hermitian form, then
$$(D,F,h,0)= {\rm Cay}(D,F,-N)$$
where $N(u)=h(u,u)$ for all $u\in F$.

 Let $S$ be a quadratic  \'etale $R$-algebra with canonical involution $\can$.
Let $(P,h)$ be a ternary nondegenerate $\can$-hermitian space such that $\bigwedge^3 (P,h)
\cong \langle 1\rangle $. Choose an isomorphism $\alpha: \bigwedge^3 (P,h)\to \langle 1\rangle $.
Define a cross product $\times_{\alpha}:P\times P\to P$ as in [T, p.~5122]. Note that
$\times_{\alpha}$ can be derived solely from $T,P,h$ and $\alpha$.
The $R$-module
${\rm Cay}(S,P,h,\alpha)=S\oplus P$ becomes an octonion algebra under the multiplication
$$(a,u)(b,v)=(ab-h(v,u), va+u\bar b+u\times_{\alpha}v)$$
for $u,v\in P$ and $a,b\in S$, with norm
$$n((a,u))=n_S(a)+h(u,u).$$
If the ternary hermitian space $(P,h)$
is orthogonally decomposable this construction is independent of
the choice of the isomorphism $\alpha$.
Any octonion algebra which contains a quadratic  \'etale  subalgebra
can be constructed like this (Petersson-Racine [PR, 3.8] or Thakur [T]).
 We obtain $$(S,P,h,\times_{\alpha})={\rm Cay}(S,P,h,\alpha).$$

 Let $2\in R^\times$. There exists a general construction for quaternion algebras over $R$ [Pu]:
Let $F$ be a finitely generated projective $R$-module of constant
rank 3 such that $\det\,F=\bigwedge^3 F\cong R$,
and $N \colon F\to R$ a nondegenerate quadratic form,
i.e., a selfdual isomorphism $\varphi_N \colon F \tilde{\to}{\check{F}},\,x\to N(x,\_)$.
(Any selfdual isomorphism $\varphi$ from $F$ to its dual
bundle ${\check{ F}}$ canonically induces a nondegenerate quadratic form on
$F$ via $N_\varphi (x) \colon ={1 \over 2} \langle x,\varphi (x) \rangle$
where $\langle \ ,\ \rangle \colon F\times {\check{ F}}
\to R$ is the canonical pairing, cf. [P2, 3.2].) Assume that $N$ has trivial determinant.
Choose an isomorphism $\alpha \colon {\det}\,{F} \tilde{\to} R$,
then $\alpha$ is unique up to an invertible factor in $R $ and induces a bilinear map
$$F\times F \to {\check{ F}}$$
$$(u,v) \to u \times_\alpha v \colon =\alpha (u \wedge v \wedge \_ )$$
as well as an isomorphism $\beta \colon \det\,{\check{F}} \tilde{\to} R$ given by the condition
$$\alpha (u_1 \wedge u_2 \wedge u )\beta (\check{u}_1 \wedge \check{u}_2
 \wedge \check{u}_3 )= \det ( \langle u_i ,\check{u}_j \rangle )$$
for $u_i$ in $F,\ \check{u}_j$ in $\check{ F},\ 1 \leq
i,j \leq 3$, where again $\langle \ ,\ \rangle \colon F\times {\check{ F}}
\to R$ is the canonical pairing. The $R$-module
$${\rm Quat}(F,N)  =R \oplus F$$
 becomes a quaternion algebra over $R$ under the multiplication
$$(a,u)(b,v)=(ab-{1\over{2}} N(u,v),\;
av+bu+{1\over{2}}\varphi_N^{-1}(u\times_\alpha v))$$
 for all $a,b$ in $R$, $u,v$ in $F$, with norm
$$N_{{\rm Quat}(F,N)}((a,u))=a^2+N (u).$$
Every quaternion algebra over $R$ is of the type  ${\rm Quat}(F,N)$.
We obtain
$${\rm Quat}(F,N)=(R,F,B,\times_\alpha),$$
where $B$ is the nondegenerate symmetric bilinear form on $F$ associated with $N$.

We have  thus unified the three constructions known for composition algebras.

\ignore{\smallskip
{\it Acknowledgements:} The author would like to acknowledge the support of the Deutsche Forschungsgemeinschaft
(in particular, the ``Georg-Thieme-Ged\"{a}chtnis--stiftung'') during her stay at the University of Trento,
 and thanks the Department of Mathematics at Trento for its hospitality and congenial atmosphere.
}

\end{document}